\documentclass[leqno]{article}
\usepackage{latexsym,amssymb,amsmath}
\usepackage[all]{xy}
\newtheorem{itheorem}{Theorem}
\newtheorem{theorem}{Theorem}[section]
\newtheorem{proposition}[theorem]{Proposition}
\newtheorem{lemma}[theorem]{Lemma}
\def\demo #1{
    \vskip-\lastskip \vskip 12pt plus 3pt minus 9pt 
\noindent{\bf #1}\enspace}
\def\noproof{{\unskip\nobreak\hfill\penalty50\hskip2em\hbox{}%
     \nobreak\hfill$\square$\parfillskip=0pt%
     \finalhyphendemerits=0\par}}
\def\enddemo{\ifmmode\eqno\square\else\noproof \vskip 12pt plus 3pt minus 9pt\fi}
\def\diagram{\renewcommand\arraystretch{1.5} $$ \begin{array}}
\def\enddiagram{\end{array} $$ \renewcommand\arraystretch{1}}

\def\rmnewname#1{\expandafter\gdef\csname#1\endcsname{{\mathop{\rm
#1}\nolimits}}}
\def\itnewname#1{{\expandafter\gdef\csname#1\endcsname{{\mathop{\it
#1}\nolimits}}}}

 \itnewname{ab}
 \itnewname{can}
 \rmnewname{codim}
 \rmnewname{cd}
 \rmnewname{cone}
 \rmnewname{coker}
 \itnewname{cor}
 \rmnewname{CH}
 \rmnewname{Cl}
 \rmnewname{div}
 \rmnewname{Div}
 \itnewname{et}
 \rmnewname{Ext}
 \rmnewname{Frob}
 \rmnewname{Hom}
 \itnewname{id}
 \rmnewname{im}
 \itnewname{inf}
 \rmnewname{Ind}
 \rmnewname{ker}
 \rmnewname{length}
 \itnewname{Mor}
 \rmnewname{Pic}
 \itnewname{pr}
 \rmnewname{Quot}
 \rmnewname{rank}
 \itnewname{real}
 \itnewname{rec}
 \itnewname{res}
 \rmnewname{rk}
 \rmnewname{Sch}
 \rmnewname{Schemes}
 \itnewname{Sign}
 \rmnewname{Sm}
 \rmnewname{SmCor}
 \rmnewname{Sp}
 \itnewname{sp}
 \rmnewname{Spec}
 \rmnewname{supp}
 \itnewname{tame}
 \rmnewname{top}
 \rmnewname{Tor}
 \rmnewname{tor}
 \rmnewname{vcd}
 \itnewname{Zar}
 
 \def\C{{\mathbb C}}
 \def\F{{\mathbb F}}

 \def\p{{\mathfrak p}}

 \def\Q{{\mathbb Q}}
 \def\R{{\mathbb R}}
 
 \def\Z{{\mathbb Z}}

 \def\zp{{\mathbb Z/p\mathbb Z}}

 \def\ms{\medskip}
 \def\ds{\displaystyle}
 
 \def\c{{\mathfrak c}}

 \def\Cal#1{{\cal #1}}

 \def\mf#1{{\mathfrak #1}}

 \def\lang{\longrightarrow}
 
 \def\mapr#1{\mathrel{\mbox{$\stackrel{#1}{\longrightarrow}$}}}
 
 \def\mapd#1{\Big\downarrow\rlap{$\vcenter{\hbox{$\scriptstyle{{#1}}$}}$}}
 
 \def\liso{\mathrel{\hbox{$\longrightarrow$} \kern-14pt\lower-4pt\hbox{$\scriptstyle\sim$}\kern7pt}}
 \def\r@iso{\mathrel{\lower2pt\hbox{$\scriptstyle\sim$} \kern-8pt\hbox{$\rightarrow$}}}
 \def\riso#1{\mathrel{\stackrel{\!\!#1}{\r@iso}}}
 \def\lr@iso{\mathrel{\kern6pt\lower2pt\hbox{$\scriptstyle\sim$} \kern-12pt\hbox{$\longrightarrow$}}}
 \def\lriso#1{\mathrel{\mathop{\lr@iso}\limits^{#1}}}

 \def\mapd#1{\Big\downarrow\rlap {$\vcenter{\hbox{$\scriptstyle{{#1}}$}}$}}
 \def\surjr#1{\stackrel{#1}{\hbox{$\relbar \! \! \! \twoheadrightarrow$}}}

 \def\surjd#1{\lower4pt\hbox{$\downarrow$}\kern-6.95pt\Big\downarrow\rlap {$\vcenter{\hbox{$\scriptstyle{{#1}}$}}$}}
 \def\surju#1{\lower-4pt\hbox{$\uparrow$}\kern-6.95pt\Big\uparrow\rlap {$\vcenter{\hbox{$\scriptstyle{{#1}}$}}$}}

 \def\dirlim{\mathop{\lim\limits_{{\longrightarrow}}\,}}
 \def\projlim{\mathop{\lim\limits_{{\longleftarrow}}\,}}

 \def\tcup{\mathop{\mathord{\cup}\mkern-8.5mu^{\hbox{.}}\mkern 5mu}}
 \def\ttcup{\mathop{\mathord{\bigcup}\mkern-10mu^{\hbox{.}}\mkern 5mu}}

 \def\hang{\hangindent\Itemindent}
 \def\textindent#1{\hskip\Itemindent\llap{\hbox{\rm #1}\enspace}\ignorespaces}
 \def\Item{\par\noindent\hang\textindent}
 
 \newdimen\Itemindent \Itemindent=.9cm
 
 \font\emar = cmsy10 scaled\magstep4 \font\emas = cmsy10 scaled\magstep2
 \def\freeproduct{\mathop{\lower.4mm\hbox{\emar \symbol{3}}}\limits}
 \def\freeproductsmall{\mathop{\hbox{$\ast$}}\limits}
 \def\freeproductmed{\mathop{\lower.2mm\hbox{\emas \symbol{3}}}\limits}
 \def\ressum{\mathop{\hbox{${\ds\bigoplus}'$}}\limits}
 
 \font\Russ=wncyr10  
 
 \def\Sha{\hbox{\Russ\char88}}

 \newcommand{\wingl}{[\![}
 \newcommand{\wingr}{]\!]}
 \newcommand{\md}{\hbox{\it \scriptsize mod}}

\title{\bf\boldmath On the relation between $2$ and $\infty$ in
 Galois cohomology of number fields}
\author{by Alexander Schmidt\\
at Heidelberg}
\date{}
\begin{document}
 \maketitle

\vskip1cm

\noindent {\footnotesize {\sc\noindent Abstract}: We remove the assumption \lq $p \not = 2$ or
$k$ is totally imaginary' from several well-known theorems on Galois groups with restricted
ramification of number fields. For example, we show that the Galois group of the maximal
extension of a number field $k$ which is unramified outside 2 has  finite cohomological
$2$-dimension (also if $k$ has real places).



}

\vskip2cm

\section{Introduction}
Number theorist's nightmare, the prime number $2$, frequently
causes technical problems and requires additional efforts. In
Galois cohomology the problems with $p=2$ are essentially due to
the fact that the decomposition groups of the real places are
$2$-groups and so the case of a totally imaginary number field is
comparatively easier to deal with.

 A classical
object of study in number theory is Galois groups with restricted
ramification. For a number field $k$, a set $S$ of primes of $k$
and a prime number $p$, one is interested in the Galois group
$G_S(p)=G(k_S(p)|k)$ of the maximal $p$-extension $k_S(p)$ of $k$
which is unramified outside $S$. If $S$ is empty, then $G_S(p)$
is the Galois group of the so-called $p$-class field tower of~$k$
and, besides the fact that it can be infinite (Golod-\v
Safarevi\v c), not much is known about this group. The situation
is easier in the case that $S$ contains the set $S_p$ of primes
dividing $p$, where the cohomological dimension of $G_S(p)$ is
known to be less than or equal to two (cf.\ \cite{NSW}, (8.3.17),
(10.4.9)). However, there is an exception: if $p=2$ and $k$ has at
least one real place. If, in this exceptional case, $S$ contains
all real places, then these places become complex in $k_S(2)$ and
therefore $G_S(2)$, containing involutions, has infinite
cohomological dimension. Furthermore, the virtual cohomological
dimension $\vcd\, G_S(2)$ is less than or equal to two in this
case, i.e.\ $G_S(2)$ has an open subgroup $U$ with $\cd\, U \leq
2$. The case when not all real places are in $S$ has been open so
far and is the subject of this paper.

\begin{itheorem} \label{haupt} Let $k$ be a number field and let $S$ be a set of
primes of\/ $k$ which contains all primes dividing $2$. If no real
prime is in $S$,  then  $\cd \, G_S(2) \leq 2$. If  $S$ contains
real primes, then they become complex in $k_S(2)$ and  $\cd \,
G_S(2)=\infty$, $\vcd\,G_S(2)\leq 2$.

If $S$ is finite, then $H^i(G_S(2)):=H^i(G_S(2),\Z/2\Z)$ is
finite for all $i$ and
\[
\chi_2(G_S(2))= -r_2,
\]
where $\chi_2(G_S(2))=\sum_{i=0}^2 (-1)^i \dim_{\F_2}
H^i(G_S(2))$ is the second partial Euler characteristic and $r_2$
is the number of complex places of\/ $k$.
\end{itheorem}

 The key for the proof of theorem~\ref{haupt} is the following
theorem~\ref{riemann} in the case $p=2$ and $T=S \cup S_\R$,
where $S_\R$ is the set of real places of $k$.
Theorem~\ref{riemann} is the number theoretical analogue of
Riemann's existence theorem and  was previously known under the
assumption that $p$ is odd or that $S$ contains $S_\R$ (see
\cite{NSW}, (10.5.1)).

\begin{itheorem}\label{riemann}
Let $k$ be a number field, $p$ a prime number and $T \supset S
\supseteq S_p$ sets of primes of $k$. Then the canonical
homomorphism
\[
\freeproductmed_{{\mathfrak p}\,\in \,T \smallsetminus S (k_S(p))}
T(k_\p(p)|k_\p) \mapr{} G(k_T(p)|k_S(p))
\]
is an isomorphism. Here $T(k_\p(p)|k_\p) \subset G(k_\p(p)|k_\p)$
is the inertia group and $\freeproductmed$ denotes the free
pro-$p$-product.
\end{itheorem}

Since the cyclotomic $\Z_2$-extension $k_\infty(2)$ of $k$ is
contained in $k_{S_2}(2)$, the group $G_{S_2}(2)$ is infinite, in
particular, it is nontrivial. Hence, for $S \supseteq S_2$ and $S
\cap S_\R=\varnothing$, the group $G_S(2)$ is of cohomological
dimension $1$ or $2$. The next theorem gives a criterion for
which case occurs. In condition $(3)$ below, $\Cl^0_S(k)(2)$
denotes the $2$-torsion part of the $S$-ideal class group in the
narrow sense of\/ $k$.

\begin{itheorem}\label{freecrit}
Assume that $S \supseteq S_2$ and $S \cap S_\R=\varnothing$. Then
$\cd\, G_S(2)=1$ if and only if the following conditions
$(1)$--$(3)$ hold. \smallskip

{ \Item{$(1)$} $S_2=\{ \p_0 \}$, i.e.\ there exist exactly one
prime dividing $2$ in $k$. \Item{$(2)$} $S= \{ \p_0\} \cup
\{\hbox{\rm complex places}\}$. \Item{$(3)$} $\Cl_S^0(k)(2)=0$.

}

\smallskip\noindent In this case, $G_S(2)$ is a free
pro-$2$-group of rank $r_2+1$ and $\p_0$ does not split in $k_{S
\cup S_\R}(2)$. In particular, if\/ $k$ is totally real and
$G_S(2)$ is free, then $k_S(2)=k_\infty(2)$.
\end{itheorem}

Let $k$ be a number field, $p$ a prime number  and $S\supseteq
S_p$ a set of places of $k$. A  (necessarily infinite) extension
$K|k$ is called $p$-$S$-closed if it has no $p$-extension which
is unramified outside $S$. If $p$ is odd and $K$ is
$p$-$S$-closed, then the group
$\Cl_{S}(K(\mu_p))(p)(j)^{G(K(\mu_p)|K)}$ is trivial for $j=0,-1$,
where $\mu_p$ is the group of $p$-th roots of unity, $(p)$
denotes the $p$-torsion part and $(j)$  the $j$-th Tate-twist
(see \cite{NSW}, (10.4.7)). The corresponding result for $p=2$ is
the following

\begin{itheorem}\label{classnumb}
Let $k$ be a number field, $S\supseteq S_2$ a set of primes of\/
$k$ and $K$ a $2$-$S$-closed extension of\/ $k$. Then the
following holds.
\begin{description}
\item{\rm (i)} $\Cl_{S} (K(\mu_4))(2)=0$.
\item{\rm (ii)} $\Cl^0_S (K) (2)=0.$
\end{description}
\end{itheorem}

\noindent {\bf Remarks:} 1. The triviality of $\Cl (K)(2)$, and
hence also that of $\Cl_S(K)(2)$, follows easily from the
principal ideal theorem;  assertions (i) and (ii) do not.

\smallskip\noindent
2. In (i) one can replace $K(\mu_4)$ by any totally imaginary
extension of degree~$2$ of $K$ in $K_S(2)$.

\bigskip

Finally, we consider the full extension $k_S$, i.e.\ the maximal
extension of $k$ which is unramified outside $S$, and its Galois
group $G_S=G(k_S|k)$.

\begin{itheorem}\label{prof}
Let $k$ be a number field and $S$ a set of primes of $k$
containing all primes dividing $2$. Then $\vcd_2 G_S \leq 2$ and
$\cd_2G_S\leq 2$ if and only if $S$ contains no real primes. For
every discrete $G_S(2)$-module $A$ the inflation maps
\[
\inf: H^i(G_S(2),A) \lang H^i(G_S,A)(2)
\]
are isomorphisms for all $i\geq 1$.
\end{itheorem}

\noindent {\bf Remark:} If $\cd\, G_S(K)(2)=2$ (e.g.\ if $K$
contains at least two primes dividing~$2$) for some finite
subextension $K$ of $k$ in $k_S$, then $\vcd_2 G_S=2$. This is
always the case if $S \supset S_\R$ because the class numbers of
the cyclotomic fields $\Q(\mu_{2^n})$ are nontrivial for $n \gg
0$. But, for example, we do not know whether $\cd_2
G(\Q_{S_2}|\Q)$ equals $1$ or $2$. The answer would be
\lq$\hspace{.5pt}2\hspace{.5pt}$\rq\ if at least one of the real
cyclotomic fields $\Q(\mu_{2^n})^+$, $n=2,3,\ldots$, would have a
nontrivial class number. But this is unknown.

\medskip In section \ref{4} we investigate the relation between the cohomology of the
group $G_S(k)$ and the modified \'{e}tale cohomology of the scheme $\Spec (\Cal O_{k,S})$. A
discrete $G_S(k)$-module $A$ induces a locally constant sheaf on $\Spec (\Cal
O_{k,S})_{\et,\md}$, which we will denote by the same letter. We show the following theorem
which is well-known if $S$ contains all real primes (and also for odd $p$).

\begin{itheorem}\label{et}
Let $k$ be a number field and $S$ a finite set of primes of\/ $k$ containing all primes
dividing $2$. Then for every $2$-primary discrete $G_S(k)$-module $A$ the natural comparison
maps
\[
H^i(G_S(k),A) \lang H^i_{\et,\md}(\Spec (\Cal O_{k,S}),A)
\]
are isomorphisms for all $i\geq 0$.
\end{itheorem}

For finite $A$ it is not difficult to show that the modified \'{e}tale cohomology groups on the
right hand side of the comparison map are finite and that they vanish for $i \geq 3$\/ if $S$
contains no real primes. Therefore one could deduce theorem~\ref{haupt} (with $G_S(k)(2)$
replaced by $G_S(k)$) from theorem~\ref{et}. However, in order to prove theorem \ref{et},
one needs information on the interaction between the decomposition groups of the real primes
and so theorem~\ref{haupt} and theorem~\ref{et} are both consequences of
theorem~\ref{riemann}.

\bigskip The main ingredients in the proofs of theorems 1--5 are
Poitou-Tate duality, the validity of the weak Leopoldt-conjecture
for the cyclotomic $\Z_p$-extension and, most essential,  the
systematic use of free products of bund\-les of profinite groups
over a topological base.  The reason that the above theorems had
not been proven earlier seems to be a psychological one. At least
the author always thought that one has to prove
theorem~\ref{haupt} first, before showing the other assertions.
For example, theorem~\ref{riemann} for $p=2$, $T=S_2 \cup S_\R$
and $S=S_2$ was known if $k_{S_2}(2)=k_\infty(2)$ (see
\cite{saf}, \S4.2 for the case $k=\Q$  and \cite{W2}, Satz 1.4
for the general case). But now it is theorem~\ref{riemann} which
is used in the proof of theorem~\ref{haupt}. Finally, we should
mention that theorem~\ref{haupt} was formulated as a conjecture
in O.\ Neumann's article \cite{N}.

\bigskip
The author wants to thank K. Wingberg for his comments which led
to a substantial simplification in the proof of
theorem~\ref{riemann}.

\section{Free products of inertia groups} \label{secfree}
In this section we briefly collect some facts on free products of
profinite groups and how they naturally occur in number theory.
For a more detailed presentation and for proofs of the facts cited
below we refer the reader to \cite{NSW}, chap.\,IV and
chap.\,X,\,\S1.

A profinite space is a topological space which is compact and
totally disconnected. Equivalently, a profinite space is a
topological inverse limit of finite discrete spaces. A profinite
group is a group object in the category of profinite spaces. It
can be shown that a profinite group is the inverse limit of
finite groups. A full class of finite groups $\c$ is a full
subcategory of the category of all finite groups which is closed
under taking subgroups, quotients and extensions. A
pro-$\c$-group is a profinite group which is the inverse limit of
groups in $\c$.

Let $T$ be a profinite space. A {\em bundle of profinite
groups}\/ $\cal G$ over $T$ is a group object in the category of
profinite spaces over $T$. We say that $\cal G$ is a bundle of
pro-$\c$-groups if the fibre ${\cal G}_t$ of $\cal G$ over every
point $t \in T$ is a pro-$\c$-group. The functor ``constant
bundle'', which assigns to a pro-$\c$-group $G$ the bundle
$\pr_2: G \times T \to T$ has a left adjoint
\[
\begin{array}{ccc}
\left\{\hbox{bundles of pro-$\c$-groups over $T$}\right\}
&\mapr{}&
\left\{\hbox{pro-$\c$-groups}\right\}\\
{\cal G}& \longmapsto & \freeproductmed_T {\cal G}.
\end{array}
\]
The image ${\freeproductmed}_T {\cal G}$ of a bundle $\cal G$
under this functor is called its free pro-$\c$-product. It
satisfies a universal property which is determined by the functor
adjunction. Bundles of pro-$\c$-groups often arise in the
following way:

 Let $G$ be a pro-$\c$-group and assume we are given a continuous family
of closed subgroups of $G$, i.e.\ a family of closed subgroups
$\{ G_t\}_{t \in T}$ indexed by the points of a profinite space
$T$ which has the property that for every open subgroup $U
\subset G$ the set $T(U)=\{t \in T\,|\, G_t \subseteq U\}$ is open
in $T$. Then
\[
{\cal G} = \{(g,t)\in G \times T\,|\, g \in G_t \}
\]
is in a natural way a bundle of pro-$\c$-groups over $T$. We have
a canonical homomorphism
\[
\phi: \freeproductmed_T {\cal G} \mapr{} G
\]
and we say that $G$ is the free product of the family $\{ G_t\}_{t \in T}$ if $\phi$ is an
isomorphism.

The usual free pro-$\c$-product of a discrete family of
pro-$\c$-groups as defined in various places in the literature
(e.g.\ \cite{Nk}) fits into the picture as follows. For a family
$\{G_i\}_{i\in I}$ we consider the disjoint union $(\tcup_i G_i)
\tcup \{\ast\}$ of the $G_i$ and one external point $\ast$.
Equipped with a suitable topology, this is a bundle of
pro-$\c$-groups over the one-point compactification $\bar I = I
\tcup \{\ast\}$ of $I$ and the free pro-$\c$-product of the family
$\{G_i\}_{i \in I}$ coincides with that of the bundle (cf.\
\cite{NSW}, chap.IV,\,\S3, examples 2 and 4). For the free product
of a discrete family of pro-$\c$-groups we have the following
profinite version of Kurosh's subgroup theorem (see \cite{NBW} or
\cite{NSW}, (4.2.1)).

\begin{theorem} \label{kurosh} Let $G= \freeproductmed_{i \in I} G_i $ be the free
pro-$\c$-product of the discrete family $G_i$ and let $H$ be an open subgroup of $G$. Then
there exist systems $S_i$ of representatives $s_i$ of the double coset decomposition $G=
\ttcup_{s_i \in S_i} H s_i G_i$ for all $i$ and a free pro-$\c$-group $F \subseteq G$ of
finite rank
\[
\rk (F) = \sum_{i \in I} [(G:H)- \# S_i] -(G:H) +1,
\]
such that the natural inclusions induce a free product
decomposition
\[
H=\freeproductmed_{i,s_i} (G_i^{s_i} \cap H) \freeproductsmall F,
\]
where $G_i^{s_i}$ $(= s_i G_i s_i^{-1})$ denotes the conjugate
subgroup.
\end{theorem}

\ms In number theory, continuous families of pro-$\c$-groups occur in the following way. For
a number field $k$ we denote the one-point compactification of the set of all places of\/
$k$ by $\Sp(k)$. The compactifying point will be denoted by $\eta_k$ and should be thought
as the generic point of the scheme $\Spec(\Cal O_k)$ in the sense of algebraic geometry or
as the trivial valuation of $k$ from the point of view of valuation theory. For an infinite
extension $K|k$, we set
\[
\Sp (K) = \projlim_{k'} \Sp(k'),
\]
where $k'$ runs through all finite subextensions of $k$ in $K$.
The complement of the (closed and open) subset of all archimedean
places of $K$ in $\Sp(K)$ is naturally isomorphic to $\Spec(\Cal
O_K )$ endowed with the constructible topology (see \cite{ega},
chap.I,\,\S7, (7.2.11) for the definition of the constructible
topology of a scheme). Let $S$ be a set of primes of $k$ and
$\bar S$ its closure in $\Sp(k)$ ($\bar S= S$ if $S$ is finite,
$\bar S= S \cup \{\eta_k\}$ if $S$ is infinite). The pre-image
$\bar S(K)$ of $\bar S$ under the natural projection $\Sp(K) \to
\Sp(k)$ is the closure of the set $S(K)$ of all prolongations of
primes in $S$ to $K$ in $\Sp(K)$.

Now assume that $M \supset K \supset k$ are possibly infinite
extensions of $k$ such that $M|K$ is Galois and $G(M|K)$ is a
pro-$\c$-group. The natural projection $\bar S (M) \to \bar S
(K)$ has a section (in fact, there are many of them). For a fixed
section $s: \bar S (K) \to \bar S (M)$ we consider the family of
inertia groups $\{ T_{s(\p)}(M|K)\}_{\p \in \bar S(K)}$, where by
convention $T_{\eta_M}=\{1\}$. Since a finite extension of number
fields is ramified only at finitely many primes, this is a
continuous family of subgroups of $G(M|K)$ indexed by $\bar
S(K)$. We obtain a natural homomorphism
\[
\phi: \freeproductmed_{\bar S(K)} T_{s(\p)}(M|K) \mapr{} G(M|K),
\]
which we also write in the form
\[
\phi:\freeproductmed_{\p \in S(K)} T_{\p}(M|K) \mapr{} G(M|K).
\]
The cohomology groups of the free product on the left hand side
with coefficients in a trivial module do not depend on the
particularly chosen section~$s$. The question, however, whether
the homomorphism $\phi$ is an isomorphism {\it does} depend
on~$s$. Moreover, if $s$ is a section for which $\phi$ is an
isomorphism, we always find a section $s'$ for which it is not,
at least if $\c$ is not the class of $p$-groups, where $p$ is a
prime number. In the case of pro-$p$-groups this pathology does
not occur because of the following easy and well-known
\begin{lemma}\label{isolemma}
Let $p$ be a prime number and let $\phi: G' \lang G$ be a
(continuous) homomorphism of pro-$p$-groups. Let $A$ be $\zp$\/
or\/ $\Q_p/\Z_p$ with trivial action. Then $\phi$ is an
isomorphism if and only if the induced homomorphism
\[
H^i(\phi,A): H^i(G,A) \lang H^i(G',A)
\]
is an isomorphism for $i=1$ and injective for $i=2$.
\end{lemma}

 In the number theoretical situation above, we have
the following formula for the cohomology of the free product with
values in a torsion group $A$ (considered as a module with
trivial action)  and for $i \geq 1$:
\[
H^i\big(\freeproductmed_{\p \in S (K)} T_\p(M|K), A\big) =
\dirlim_{k'} \bigoplus_{\p \in S(k')} H^i(T_\p(M'|k'), A),
\]
where $k'$ runs through all finite subextensions of $k$ in $K$
and $M'$ is the maximal pro-$\c$ Galois subextension of $M|k'$
(so $M = \dirlim M'$). The limit on the right hand side depends
on $K$ and not on $k$ and we denote it by
\[
{\ressum_{\p \in S(K)}} H^i(T_\p (M|K),A).
\]
If $K|k$ is Galois, then this limit is the maximal discrete
$G(K|k)$-submodule of the product ${\prod_{\p \in S(K)}} H^i(T_\p
(M|K),A)$.

\pagebreak
\section{Proof of theorem~\ref{riemann}}

Let us first remark that for ${\mathfrak p}\,\in \,T
\smallsetminus S (k)$ the inertia group has the following
structure:
 \Item{-} if $\p$ is nonarchimedean and $N(\p)\equiv 1 \bmod p$
(i.e.\ if there is a primitive $p$-th root of unity in $k_\p$),
then $T(k_\p(p)|k_\p)$ is a free pro-$p$-group of rank~$1$, i.e.\
isomorphic to $\Z_p$.
 \Item{-} if $\p$ is nonarchimedean and $N(\p)\not \equiv 1 \bmod p$,
 then $T(k_\p(p)|k_\p)=\{ 1\}$.
 \Item{-} if $\p$ is real and $p=2$, then $T(k_\p(p)|k_\p)\cong \Z/2\Z$.
 \Item{-} if $\p$ is real and $p\not =2$ or if $\p$ is complex, then $T(k_\p(p)|k_\p)= \{ 1\}$.

\bigskip
If $p$ is odd or if $p=2$ and $S \supset S_\R$, then
theorem~\ref{riemann} is known (see \cite{NSW}, (10.5.1)). So we
assume that $p=2$ and $S \not \supset S_\R$. For a pro-$2$-group
$G$ we use the notation $H^i(G)$ for $H^i(G,\Z/2\Z)$. We start
with the following

\begin{lemma} \label{isocrit2}
Let $G$ and $G'$ be pro-$2$-groups which are generated by
involutions and assume that
$H^2(G,\Q_2/\Z_2)=0=H^2(G',\Q_2/\Z_2)$. Let $\phi: G' \to G$ be a
(continuous) homomorphism. Then the following assertions are
equivalent.
\begin{description}
\item{\rm (i)} $\phi$ is an isomorphism.
\item{\rm (ii)} $H^1(\phi): H^1(G) \to H^1(G')$ is an
isomorphism.
\item{\rm (iii)} $H^2(\phi): H^2(G) \to H^2(G')$ is an
isomorphism.
\end{description}
\end{lemma}

\demo{Proof:} Clearly, (i) implies (ii) and (iii) and, by lemma
\ref{isolemma}, (ii) and (iii) together imply (i). So it remains
to show that (ii) and (iii) are equivalent. Since
$H^2(G,\Q_2/\Z_2)=0$, the exact sequence $0 \to \Z/2\Z \to
\Q_2/\Z_2 \to \Q_2/\Z_2 \to 0$ induces the four term exact
sequence
\[
0 \to H^1 (G) \stackrel{\alpha}{\to} H^1(G, \Q_2/\Z_2)
\stackrel{\beta}{\to} H^1(G,\Q_2/\Z_2) \stackrel{\gamma}{\to}
H^2(G) \to 0.
\]
Since $G$ is generated by involutions, $\alpha$ is an
isomorphism. Hence $\beta$ is zero and $\gamma$ is an
isomorphism. The same argument also applies to $G'$ and therefore
(ii) and (iii) are both equivalent to
\begin{description}
\item{\rm (iv)} $H^1(\phi, \Q_2/\Z_2): H^1(G,\Q_2/\Z_2) \to H^1(G',\Q_2/\Z_2)$ is an
isomorphism.
\end{description}
This concludes the proof.\enddemo

We show theorem \ref{riemann} first in the special case $T=S_2
\cup S_\R$, $S=S_2$. The groups ${\freeproductmed}_{\p \in
S_\R(k_{S_2}(2))} T(k_\p(2)|k_\p)$ and $G(k_{S_2 \cup
S_\R}(2)|k_{S_2}(2))$ are both generated by involutions. Since
$H^2(T(k_\p(2)|k_\p),\Q_2/\Z_2)=0$ for every $\p \in
S_\R(k_{S_2}(2))$, we have
\[
H^2(\freeproductmed_{\p \in S_\R(k_{S_2}(2))}
T(k_\p(2)|k_\p),\Q_2/\Z_2)=0.
\]
By \cite{NSW}, (10.4.8), the inflation map
\[
H^2\big(G(k_{S_2 \cup S_\R}(2)|k_{S_2}(2)),\Q_2/\Z_2\big) \lang
H^2\big(G(k_{S_2 \cup S_\R}|k_{S_2}(2)),\Q_2/\Z_2\big)
\]
is an isomorphism and, since $k_{S_2}(2)$ contains the cyclotomic
$\Z_2$-extension $k_\infty(2)$ of $k$, the validity of the weak
Leopoldt-conjecture for the cyclotomic $\Z_p$-extension (see
\cite{NSW}, (10.3.25)) implies (by \cite{NSW},  (10.3.22)) that
\[
H^2(G(k_{S_2 \cup S_\R}(2)|k_{S_2}(2)),\Q_2/\Z_2)=0.
\]
By lemma \ref{isocrit2} and the calculation of the cohomology of
free products (see \S1), it therefore suffices to show that the
natural map
\[
H^2(\phi):  H^2\big(G(k_{S_2 \cup S_\R}(2)|k_{S_2}(2)\big)\to
\ressum_{\p \in S_\R(k_{S_2}(2))} H^2\big(T(k_\p(2)|k_\p)\big)
\]
is an isomorphism. Now let $K$ be a finite extension of $k$ inside $k_S(2)$. The $9$-term
exact sequence of Poitou-Tate induces the exact sequence \pagebreak[3]
\[
0 \to \Sha^2(K_{S_2 \cup S_\R},\Z/2\Z) \to H^2(G(k_{S_2 \cup S_\R}
|K),\Z/2\Z) \to
\]
\begin{flushright}
$\ds \bigoplus_{\p \in S_2 \cup S_\R(K)} H^2(G(\bar
k_\p|K_\p),\Z/2\Z) \to H^0(G(k_{S_2 \cup S_\R}|K),\mu_2)^\vee\to
0$,
\end{flushright}
where $\scriptstyle \vee$ denotes the Pontryagin dual.
Furthermore, we have
\[\Sha^2(K_{S_2 \cup S_\R},\Z/2\Z)\cong\Sha^1(K_{S_2 \cup
S_\R},\mu_2)^\vee\!=\! \Sha^1(K_{S_2 \cup S_\R},\Z/2\Z)^\vee
\!=\!\Cl_{S_2}(K)/2.
\]
For a finite, nontrivial extension $K'$ of $K$ inside $k_{S_2}(2)$
the corresponding homomorphism $H^0(G(k_{S_2 \cup
S_\R}|K),\mu_2)^\vee \to H^0(G(k_{S_2 \cup S_\R}|K'),\mu_2)^\vee $
is the dual of the norm map, hence trivial. Furthermore,
$H^2\big(G(\bar k_\p|(k_{S_2}(2))_\p),\Z/2\Z\big)=0$ for $\p \in
S_2(k_{S_2}(2))$ (see \cite{NSW}, (7.1.8)(i)). Therefore we
obtain the following exact sequence in the limit over all finite
subextensions $K|k$ in $k_{S_2}(2)|k$ (the omitted coefficients
are $\Z/2\Z$):
\[
\Cl_{S_2}(k_{S_2}(2))/2 \hookrightarrow H^2\big(G(k_{S_2 \cup
S_\R}|k_{S_2}(2))\big)\twoheadrightarrow\!\!\!\!\!\!\!
\ressum_{\p \in S_\R(k_{S_2}(2))}\!\!\!\!\!\! H^2\big(G(\bar
k_\p|k_\p)\big).
\]
 The principal ideal
theorem implies that $\Cl(k_{S_2}(2))(2)=0$, and therefore also
$\Cl_{S_2}(k_{S_2}(2))/2=0$. Furthermore, $G(\bar
k_\p|k_\p)=T(k_\p(2)|k_\p)$ for $\p\in S_\R(k_{S_2}(2))$ and the
inflation map
\[
H^2\big(G(k_{S_2 \cup S_\R}(2)|k_{S_2}(2))\big) \lang
H^2\big(G(k_{S_2 \cup S_\R}|k_{S_2}(2))\big)
\]
is an isomorphism (see  \cite{NSW}, (10.4.8)). This concludes the
proof of theorem~\ref{riemann} in the case $T=S_2 \cup S_\R$,
$S=S_2$. For the proof in the general case we need the

\begin{proposition} \label{freeprodext} Let $k$ be a number field, $p$ a prime number
and $T \supset S \supseteq S_p$ sets of primes in $k$. Let $K$ be
a $p$-$S_p$-closed extension of\/ $k$. Then the following
assertions are equivalent.
\begin{description}
\item{\rm (i)} The natural homomorphism
\[
\phi_{T,S_p}:\freeproductmed_{\p \in T \smallsetminus S_p(K)}
T(K_\p(p)|K_\p) \to G(K_T(p)|K)
\]
is an isomorphism.
\item{\rm (ii)} The natural homomorphisms
\[
\phi_{T,S}:\freeproductmed_{\p \in T \smallsetminus S (K_S(p))}
T(K_\p(p)|K_\p) \to G(K_T(p)|K_S(p))
\]
and
\[
\phi_{S,S_p}:\freeproductmed_{\p \in S \smallsetminus S_p (K)}
T(K_\p(p)|K_\p) \to G(K_S(p)|K)
\]
are isomorphisms.
\end{description}
Here $\freeproductmed$ denotes the free pro-$p$-product.
\end{proposition}

\demo{Proof:} If $\phi_{T,S_p}$ is an isomorphism, then also
$\phi_{S,S_p}$ is an isomorphism. Furthermore, a straightforward
application of theorem~\ref{kurosh} shows that also $\phi_{T,S}$
is an isomorphism in this case. Let us show the converse
statement. Assume that $\phi_{T,S}$ and $\phi_{S,S_p}$ are
isomorphisms. Note that all primes in $S \smallsetminus S_p
(K_S(p))$ split completely in $K_T(p)|K_S(p)$. Therefore the
extension of pro-$p$-groups
\begin{equation}\label{split}
 1 \to G(K_T(p)|K_S(p)) \to G(K_T(p)|K) \to G(K_S(p)|K) \to 1
\end{equation}
splits. By lemma~\ref{isolemma}, we have to show that the induced
homomorphism
\[
H^i(\phi_{T,S_p}): H^i\big(G(K_T(p)|K)\big) \lang \ressum_{\p \in
T\smallsetminus S_p(K)} H^i\big(T(K_\p(p)|K_\p)\big)
\]
is an isomorphism for $i=1$ and injective for $i=2$ (coefficients
$\Z/p\Z$). This follows easily from the Hochschild-Serre spectral
sequence associated to the split exact sequence (\ref{split}):
\[
E_2^{ij}= H^i\big(G(K_S(p)|K), H^j(G(K_T(p)|K_S(p)))\big)
\Longrightarrow H^{i + j}(G(K_T(p)|K)).
\]
First of all, the differentials $d_2$ are zero ($-d_2$ is the
cup-product with the extension class, see \cite{NSW}, (2.1.8)).
Furthermore,  every prime in $T\smallsetminus S(K)$ splits
completely in $K_S(p)|K$ because these primes are unramified in
$K_S(p)|K$  and $K$ contains $K_\infty (p)$. Since $\phi_{T,S}$
is an isomorphism, the $G(K_S(p)|K)$-module ($j \geq 1$)
\[
\renewcommand{\arraystretch}{1.8}
\begin{array}{ccr}
 H^j (G(K_T(p)|K_S(p)))&=& \ressum_{\p \in T\smallsetminus S (K_S(p))}
 H^j(T(K_\p(p)|K_\p))\\
 &=& \Ind_{G(K_S(p)|K)} \ressum_{\p \in T\smallsetminus S(K)}
 H^j(T(K_\p(p)|K_\p))
\end{array}
\renewcommand{\arraystretch}{1}
\]
is cohomologically trivial. Therefore we obtain short exact
sequences
\[
0 \to H^i(K_S(p)|K) \to H^i(K_T(p)|K) \to \ressum_{\p \in
T\smallsetminus S(K)}
 H^i(T(K_\p(p)|K_\p)) \to 0
\]
for $i=1,2$, and the result follows from the five-lemma.
\enddemo

Now we can prove theorem~\ref{riemann} in the general case. It is
true for odd $p$ and for $p=2$ in the special cases  $T=S_2 \cup
S_\R$, $S=S_2$ and $T=\{ \hbox{all primes}\}$, $S= S_2 \cup S_\R$.
Applying proposition~\ref{freeprodext} in the situation $p=2$,
$T=\{ \hbox{all primes}\}$, $S= S_2 \cup S_\R$ and
$K=k_{S_2}(2)$, we obtain theorem~\ref{riemann} in the \lq
extremal\rq\ case $T=\{ \hbox{all primes}\}$, $S= S_2$. Applying
proposition~\ref{freeprodext} again, we obtain the case $T=\{
\hbox{all primes}\}$ and $S$ arbitrary and then the general case.
This concludes the proof of theorem~\ref{riemann}.

\bigskip
A straightforward limit process shows the following variant of
theorem~\ref{riemann}.

\bigskip\noindent
{\bf Theorem 2' } {\it Let $k$ be a number field, $p$ a prime
number and $T \supset S \supseteq S_p$ sets of primes of $k$. Let
$K$ be a $p$-$S$-closed extension field of $k$. Then the canonical
homomorphism
\[
\freeproductmed_{{\mathfrak p}\,\in \,T \smallsetminus S (K)}
T(K_\p(p)|K_\p) \mapr{} G(K_T(p)|K)
\]
is an isomorphism. }

\section{Proofs of the remaining statements} \label{remain}
In order to prove theorem \ref{haupt}, we may assume that $S\not
\supset S_\R$ and we investigate the Hochschild-Serre spectral
sequence
\[
E_2^{ij}= H^i\big(G_{S}(2),H^j(G(k_{S\cup
S_\R}(2)|k_{S}(2))\big)\Longrightarrow H^{i+j}(G_{S\cup S_\R}(2)),
\]
where the omitted coefficient are $\Z/2\Z=\mu_2$. By theorem
\ref{riemann}, we have complete control over the $G_S(2)$-modules
$H^j(G(k_{S\cup S_\R}(2)|k_{S}(2))$, which are for $j\geq 1$
isomorphic to
\[
\Ind_{G_{S}(2)} \bigoplus_{\p \in S_\R\smallsetminus S(k)}
H^j(G(\C |\R)).
\]
In particular, $E_2^{ij}=0$ for $ij\not =0$. Therefore the
spectral sequence induces an exact sequence
\begin{equation}\label{1to2}
0 \to H^1(G_S(2)) \to H^1(G_{S \cup S_\R}(2)) \to \bigoplus_{\p
\in S_\R \smallsetminus S(k)} H^1(G(\C|\R)) \to
\end{equation}
\begin{flushright}
$ H^2(G_S(2)) \to H^2(G_{S \cup S_\R}(2)) \to \ds\bigoplus_{\p \in
S_\R \smallsetminus S(k)} H^2(G(\C|\R)) \to 0$
\end{flushright}
and exact sequences
\begin{equation}\label{ab3}
 0 \to H^i(G_S(2)) \to H^i(G_{S \cup S_\R}(2)) \to \ds\bigoplus_{\p \in
S_\R \smallsetminus S(k)} H^i(G(\C|\R)) \to 0.
\end{equation}
for $i \geq 3$. If $S$ is finite, this shows the finiteness
statement on the cohomology of $G_S(2)$ and that
\[
\chi_2(G_S(2))=\chi_2(G_{S\cup S_\R}(2)).
\]
But $\chi_2(G_{S\cup S_\R}(2))=\chi_2(G_{S\cup S_\R})=-r_2$ (see
\cite{NSW}, (8.6.16) and (10.4.8)).

\medskip For arbitrary $S$ and $i\geq 3$  the restriction map
\[
H^i(G_{S \cup S_\R}(2)) \to \bigoplus_{\p \in S_\R (k)}
H^i(G(\C|\R))
\]
is an isomorphism (see \cite{NSW},  (8.6.13)(ii) and (10.4.8)).
This together with (\ref{ab3}) shows that the natural homomorphism
\[
H^i(G_{S}(2)) \to \bigoplus_{\p \in S \cap S_\R (k)} H^i(G(\C|\R))
\]
is an isomorphism for $i\geq 3$. Therefore $\cd\ G_S(2) \leq 2$ if
$S \cap S_\R =\varnothing$. For later use we formulate the last
result as a proposition.

\begin{proposition}\label{isoab3}
Let $k$ be a number field and $S \supset S_2$ a set of primes.
Then the natural homomorphism
\[
H^i(G_{S}(2),\Z/2\Z) \to \bigoplus_{\p \in S \cap S_\R (k)}
H^i(G(\C|\R),\Z/2\Z)
\]
is an isomorphism for $i\geq 3$.
\end{proposition}

In order to conclude the proof of theorem \ref{haupt}, it remains
to show that every real prime in $S$ ramifies in $k_S(2)$. Let
$S^f$ be the subset of nonarchimedean primes in $S$. Then
theorem~\ref{riemann} yields an isomorphism
\[
\freeproductmed_{\p \in S_\R(k_{S^f}(2))} T(k_\p(2)|k_\p)\cong
G(k_S(2)|k_{S^f}(2))
\]
which shows the required assertion. This finishes the proof of
theorem~\ref{haupt}.

\vskip1cm Now we prove theorem \ref{freecrit}.  To fix conventions, we recall the following
definitions. For a set $S$ of primes of $k$ the group $\Cal O_{k,S}^\times$ of $S$-units is
defined as the subgroup in $k^\times$ of those elements which are units at every finite prime
not in $S$ and positive at every real prime not in $S$. The $S$-ideal class group
$\Cl_S^0(k)$ in the narrow sense of $k$ is the quotient of the group  of fractional ideals
of $k$ by the subgroup generated by the nonarchimedean primes in $S$ and the principal
ideals $(a)$ with $a$ positive at every real place of $k$ not contained in $S$. In
particular, $\Cl^0_{\varnothing}(k)=\Cl^0(k)$ is the ideal class group in the narrow sense
and $\Cl^0_{S \cup S_\R}(k)=\Cl_S(k)$ is the usual $S$-ideal class group. By class field
theory, $\Cl^0_S(k)$ is isomorphic to the Galois group of the maximal abelian extension of
$k$ which is unramified outside $S_\R$ and in which every prime in $S$ splits completely. By
Kummer theory, we can replace condition (3) of theorem~\ref{freecrit} by the following
condition

\smallskip\noindent
(3') \hspace{.5cm} $  \{ x\in k^\times\,\big|\,x \in
k_{\p_0}^{\times 2} \hbox{ and } 2 \,|\, v_\p(x) \hbox{ for every
finite prime $\p$}\}= k^{\times 2}.$

\begin{lemma}
If $S \supseteq S_2$ and  $\cd\,G_{S_2}(2)=1$, then $S=S_2$.
\end{lemma}

\demo{Proof:} By theorem \ref{riemann}, we have an isomorphism
\[
\freeproductmed_{\p \in S \smallsetminus S_2(k_{S_2}(2))}
T(k_\p(2)|k_\p) \mapr{\sim} G(k_S(2)|k_{S_2}(2))
\]
Since for nonarchimedean primes $\p \not\in S_2$ the maximal
unramified $2$-extension of $k_\p$ is realized by $k_\infty(2)
\subset k_{S_2}(2)$, this shows that for $\p \in S \smallsetminus
S_2$ the maximal $2$-extension of the local field $k_\p$ is
realized by $k_S(2)$ or, in other words, the natural homomorphism
\[
G(k_\p(2)|k_\p) \lang G_S(2)
\]
is injective. But for these primes we have $\cd\,
G(k_\p(2)|k_\p)\geq 2$ which shows that $S \smallsetminus
S_2=\varnothing$.
\enddemo

Now assume that $G_{S_2}(2)$ is free. For a prime $\p$ we denote
the local group $G(k_\p(2)|k_\p)$ by $\Cal G_\p$ and the inertia
group $T(k_\p(2)|k_\p)$ by $\Cal T_\p$. By \v Cebotarev's density
theorem, we find a finite set of nonarchimedean primes $T\supset
S_2$ such that the natural homomorphism
\[
H^1(G_{S_2}) \lang \bigoplus_{\p \in T \smallsetminus S} H^1(\Cal
G_\p / \Cal T_\p)
\]
is an isomorphism. It is then an easy exercise using lemma
\ref{isolemma} to show that the natural homomorphism
\[
\freeproductmed_{\p \in T \smallsetminus S_2} \Cal G_\p /\Cal T_\p
\lang G_{S_2}(2)
\]
is an isomorphism. Theorem \ref{riemann} for $T=S_2 \cup S_\R$
and $S=S_2$ and the same arguments as in the proof of proposition
\ref{freeprodext}  show that the natural homomorphism
\[
\freeproductmed_{\p \in T \smallsetminus S_2} \Cal G_\p /\Cal T_\p
\;\; \freeproductsmall \freeproductmed_{\p \in S_\R} \Cal G_\p
\lang G_{S_2 \cup S_\R}(2)
\]
is an isomorphism. Then, by (\cite{W}, Theorem 6) or (\cite{NSW},
 (10.7.2)), we obtain the conditions (1)--(3) and that the
unique prime $\p_0$ dividing $2$ in $k$ does not split in $k_{S_2
\cup S_\R}$. If, on the other hand, conditions (1)--(3) of
theorem~\ref{freecrit} are satisfied, then we obtain (loc.\ cit.)
the above isomorphism and  deduce that $G_{S_2}(2)$ is free. The
statement on the rank of $G_{S_2}(2)$  follows from
$\chi_2(G_{S_2}(2))= -r_2$. If $k$ is totally real, then the
homomorphism
\[
G_{S_2}(2) \lang G(k_\infty(2)|k)
\]
is a surjection of free pro-$2$-groups of rank $1$ and hence an
isomorphism. This concludes the proof of theorem~\ref{freecrit}.

\vskip1cm Next we show theorem~\ref{classnumb}. Let $S$ be a set
of finite primes of $k$ and $\Sigma=S \cup S_\R$. If $S$ is
finite, then the image of the group of $\Sigma$-units of $k$
under the logarithm map $\hbox{\it Log}: \Cal O_{k,\Sigma}^\times
 \lang \bigoplus_{v \in \Sigma} \R$, $a \mapsto (\log |a|_v)_{v
\in S}$ is a lattice of rank equal to $\# S + r_1 +r_2 -1$
(Dirichlet's unit theorem). Complementary to this map is the
signature map (which is also defined for infinite $S$)
\[
\Sign_{k,S}: \Cal O_{k,\Sigma}^\times   \lang \bigoplus_{v \in
S_\R} \R^\times/\R^{\times 2}.
\]
More or less by definition, there exists a five-term exact
sequence
\[ 0 \to
\Cal O_{k,S}^\times   \rightarrow \Cal O_{k,\Sigma}^\times  \to
\bigoplus_{v \in S_\R(k)} \R^\times/\R^{\times 2} \to \Cl^0_S(k)
\rightarrow \Cl^0_\Sigma (k) \to 0,
\]
and so the cokernel of $\Sign_{k,S}$ measures the difference
between the usual $S$-ideal class group $\Cl_S(k)=\Cl^0_\Sigma
(k)$ and that in the narrow sense. Of course this discussion is
void if $k$ is totally imaginary. If $K$ is an infinite extension
of $k$, we define the signature map
\[
\Sign_{K,S}: \Cal O_{K,\Sigma}^\times \lang \dirlim_{k'}
\bigoplus_{v \in S_\R (k')} \R^\times /\R^{\times}
\]
as the limit over the signature maps $\Sign_{k',S}$, where $k'$
runs through all finite subextension $k'|k$ of $K|k$. If $K$ is
$2$-$S$-closed, then $\Cl_S(K)(2)=0$ and so statement (ii) of
theorem~\ref{classnumb} is equivalent to the statement that
$\Sign_K$ is surjective.

\medskip Now assume that $k$, $S$, $K$ are as in
theorem~\ref{classnumb}. By theorem~\ref{haupt}, all real places
in $S$ become complex in $K$. By the principal ideal theorem,
$\Cl(K)(2)=2$ and so statement (i) and (ii) are trivial if $K$ is
totally imaginary (note that $K=K(\mu_4)$ in this case).  So we
may assume that $S_\R(K)\neq \varnothing$ and, by
theorem~\ref{haupt}, we may suppose $S \cap S_\R= \varnothing$.

Let $K'=K(\mu_4)$. Then $K'$ is totally imaginary and $G=G(K'|K)$
is cyclic of order~$2$. Let $\Sigma=S \cup S_\R$ and let
$K_\Sigma$ be the maximal (not just the pro-$2$) extension of $K$
which is unramified outside $\Sigma$.  Inspecting  the
Hochschild-Serre spectral sequence associated to
$K_\Sigma|K_\Sigma(2)|K$ and using the well-known calculation of
$H^i(G(K_\Sigma|K), \Cal O^\times_{K_\Sigma,\Sigma})$ (cf.\
\cite{NSW},  (10.4.8)) we see that
\begin{equation} \label{eq1}
\begin{array}{ccl}
H^1\big(G(K_\Sigma(2)|K), \Cal O_{K_\Sigma(2),\Sigma}^\times\big)&
=&H^1\big(G(K_\Sigma|K), \Cal O_{K_\Sigma,\Sigma}^\times \big)(2)\\
&=& \Cl_S(K)(2)=0
\end{array}
\end{equation}
and the same argument shows that
\begin{equation}\label{eq2}
 H^1\big(G(K_\Sigma(2)|K'), \Cal
O_{K_\Sigma(2),\Sigma}^\times)\cong \Cl_S(K')(2).
\end{equation}
 Next we consider the Hochschild-Serre spectral sequence for the
extension\linebreak $K_\Sigma(2)|K'|K$ and the module $\Cal
O^\times_{K_{\Sigma}(2),\Sigma}$. By (\ref{eq1}) and (\ref{eq2}),
we obtain an exact sequence
\[
0 \to \Cl_S(K')(2)^G \to H^2(G, \Cal
O_{K',\Sigma}^\times)\stackrel{\phi}{ \to}
H^2\big(G(K_\Sigma(2)|K), \Cal O_{K_\Sigma(2),\Sigma}^\times\big).
\]
Since $G$ is a $2$-group, in order to prove assertion (i), it
suffices to show that $\phi$ is injective. Let $c$ be a generator
of the cyclic group $H^2(G,\Z)$. For each prime $\p\in S_\R (K)$
(respectively for the chosen prolongation of $\p$ to
$K_\Sigma(2)$, cf.\ the discussion in section~1), the composition
$T_\p(K_\Sigma(2)|K) \to G(K_\Sigma(2)|K)\to G$ is an isomorphism
and we denote the image of $c$ in $H^2(T_\p(K_\Sigma(2)|K),\Z)$
by $c_\p$. As is well known, the cup-product with $c$ induces an
isomorphism $\hat H^0(G,\Cal
O_{K',\Sigma}^\times)\stackrel{\sim}{\to}H^2(G,\Cal
O_{K',\Sigma}^\times)$ and the similar statement holds for each
$c_\p$, $\p \in S_\R(K)$.

The quotient $\Cal O_{K_\Sigma(2),\Sigma}^\times / \mu_{2^\infty}$
is uniquely $2$-divisible, and so we obtain a natural isomorphism
\[
H^2\big(G(K_\Sigma(2)|K),
\mu_{2^\infty}\big)\stackrel{\sim}{\lang}
H^2\big(G(K_\Sigma(2)|K), \Cal O_{K_\Sigma(2),\Sigma}^\times\big).
\]
Furthermore, for each $\p \in S_\R \smallsetminus S$ we obtain an
isomorphism
\[
\renewcommand{\arraystretch}{1.5}
\begin{array}{ccl}
 H^2(T_\p(K_\Sigma(2)|K),\mu_{2^\infty})&\stackrel{\sim}{\to}&
H^2(T_\p(K_\Sigma(2)|K),\Cal O_{K_\Sigma(2),\Sigma}^\times)\\
& \cong & H^2(G(\bar K_\p|K_\p), \bar K^\times_\p).
\end{array}
\renewcommand{\arraystretch}{1}
\]
Therefore, the calculation of the cohomology  in dimension $i\geq
2$  of free products with values in torsion modules (see
\cite{N}, Satz 4.1 or  \cite{NSW}, (4.1.4)) and
theorem~\ref{riemann} for the pair $\Sigma$, $S$ show that we
have a natural isomorphism
\[
H^2(G(K_\Sigma(2)|K), \Cal O_{K_\Sigma(2),\Sigma}^\times)
\stackrel{\sim}{\to} \ressum_{\p \in S_\R (K)} H^2(G(\bar
K_\p|K_\p), \bar K_\p^\times).
\]
(Alternatively, we could have obtained this isomorphism from the
calculation of the cohomology of the $\Sigma$-units, cf.\
(\cite{NSW}, (8.3.10)(iii)) by passing to the limit over all
finite subextensions of $k$ in $K$). We obtain the following
commutative diagram

{\footnotesize \vskip-.7cm
\[
 \xymatrix{
 \hat H^0(G, \Cal O_{K',\Sigma}^\times) \ar@{-->}[rr]^{\psi}\ar[d]^{\cup\, c}_\wr &&
 \ressum_{\p \in S_\R (K)}^{\phantom{K}}\!\! \! \hat H^0(G(\bar K_\p|K_\p), \bar
 K_\p^\times)\ar[d]^{\ressum \cup\, c_\p}_\wr\\
 H^2(G, \Cal O_{K',\Sigma}^\times ) \ar [r]^(.4){\phi}&
 H^2(G(K_\Sigma(2)|K),\Cal O_{K_\Sigma(2),\Sigma}^\times)\ar[r]^(.5){\sim}&
 \ressum_{\p \in S_\R (K)}^{\phantom{K}}\!\!\! H^2(G(\bar K_\p|K_\p), \bar
 K_\p^\times).}
 \]
}

\noindent Hence $ \ker (\phi)\cong \ker(\psi)$ and
$\coker(\phi)\cong \coker(\psi)$. Since $\hat H^0(G, \Cal
O_{K',\Sigma}^\times )= \break\Cal O_{K,\Sigma}^\times
/N_{K'|K}(\Cal O_{K',\Sigma}^\times)$, each element in
$\ker(\psi)$ is represented by an $S$-unit in $K$ and we have to
show that all these are norms of  $\Sigma$-units in $K'$. Let
$e\in \Cal O_{K,S}^\times$. Then $K(\sqrt{e}\,)|K$ is a
$2$-extension which is unramified outside $S$, hence trivial.
Therefore $e$ is a square in $K$ and if $f^2=e$, then $f\in \Cal
O_{K,\Sigma}^\times$ and $e=N_{K'|K} (f)$. This concludes the
proof of assertion (i).

\noindent To show assertion (ii), it remains to show that $\coker
(\Sign_{K,S})=\coker (\psi) \cong \coker (\phi)$ is trivial.
Using the same spectral sequence as before, in order to see that
$\coker (\phi)=0$, it suffices to show that the spectral terms
\begin{description}
\item{-} $E_2^{02}= H^0\big(G, H^2(G(K_\Sigma(2)|K'),\Cal O_{K_\Sigma(2),\Sigma}^\times
)\big)\;$ and
\item{-} $E_2^{11}=H^1(G, \Cl_S(K')(2))$
\end{description}
are trivial. The first assertion is easy, because $K'$ is totally
imaginary and contains $k_\infty(2)$ and so
$H^2(G(K_\Sigma(2)|K'),\Cal O_{K_\Sigma(2),\Sigma}^\times)
=\nobreak 0$. That the second spectral term is trivial follows
from (i). This completes the proof of theorem~\ref{classnumb}.

\vskip1cm Finally, we prove theorem~\ref{prof}. The statement on
$\cd_2G_S$ and $\vcd_2 G_S$ follows by choosing a $2$-Sylow
subgroup $H \subset G_S$ and applying theorem~\ref{haupt} to all
finite subextensions of $k$ in $(k_S)^H$. It remains to show the
statement on the inflation map. It is equivalent to the statement
that
\[
\inf \otimes \Z_{(2)}: H^i(G_S(2),A) \otimes \Z_{(2)} \lang
H^i(G_S,A)\otimes \Z_{(2)}
\]
is an isomorphism for every discrete $G_S(2)$-module $A$ and all
$i\geq 0$, where $\Z_{(2)}$ denotes the localization of $\Z$ at
the prime ideal $(2)$.

Since cohomology commutes with inductive limits, we may assume
that $A$ is finitely generated (as a $\Z$-module). Using the
exact sequences
\[
0 \lang \tor(A) \lang A \lang A/\tor (A) \lang 0,
\]
\[
0 \lang A/\tor (A)\lang  \big( A/\tor (A)\big) \otimes \Q  \lang
\big( A/\tor (A)\big) \otimes \Q/\Z \lang 0
\]
and using the limit argument for $\big(A/\tor (A)\big) \otimes \Q/\Z $ again, we are reduced
to the case that $A$ is finite. Every finite $G_S(2)$-module is the direct sum of its
$2$-part and its prime-to-$2$-part. The statement is obvious for the prime-to-$2$-part and
every finite $2$-primary $G_S(2)$-module has a composition series whose quotients are
isomorphic to $\Z/2\Z$. Therefore we are reduced to showing the statement on the inflation
map for $A=\Z/2\Z$. But it is more convenient to work with $A=\Q_2/\Z_2$ (with trivial
action) which is possible by the exact sequence
\[
0 \lang \Z/2\Z \lang \Q_2/\Z_2 \lang \Q_2 /\Z_2 \lang 0.
\]
Using the Hochschild-Serre spectral sequence for the extensions
$k_S|k_S(2)|k$, we thus have to show that
\[
H^i(G(k_S|k_S(2)),\Q_2/\Z_2)=0
\]
for $i \geq 1$. The case $i=1$ is obvious by the definition of
the field $k_S(2)$. By theorem~\ref{haupt}, every real prime in
$S$ becomes complex in $k_S(2)$ and therefore $\cd_2
G(k_S|k_S(2))\leq 2$. It remains to show that
$H^2(G(k_S|k_S(2)),\Q_2/\Z_2)=0$. Therefore the next proposition
implies the remaining statement of theorem~\ref{prof}.

\begin{proposition}
Let $k$ be a number field, $S\supseteq S_2$ a set of primes in $k$
and $K \supseteq k_\infty(2)$ an extension of $K$ in $k_S$. Then
\[
H^2(G(k_S|K),\Q_2/\Z_2)=0.
\]
\end{proposition}

\demo{Proof:} Let $H$ be a $2$-Sylow subgroup in $G(k_S|K)$ and
$L=(k_S)^H$. Then the restriction map
\[
H^2(G(k_S|K),\Q_2/\Z_2) \lang H^2(G(k_S|L),\Q_2/\Z_2)
\]
is injective and so, replacing $K$ by $L$, we may suppose that
$k_S=K_S(2)$. Applying theorem~2' to the $2$-$S$-closed field
$K_S(2)$, we obtain an isomorphism
\[
G(K_{S \cup S_\R}(2)|K_S(2)) \cong \freeproductmed_{\p \in
S_\R(K_S(2))} T_\p(K_\p(2)|K_\p).
\]
Hence we have complete control over the Hochschild-Serre spectral
sequence associated to $K_{S \cup S_\R}(2)|K_S(2)|K$. Furthermore,
the weak Leopoldt conjecture holds for the cyclotomic
$\Z_2$-extension and $K\supseteq k_\infty (2)$, which implies that
$H^2(G(K_{S\cup S_\R}(2)|K),\Q_2/\Z_2)=0$. The exact sequence
(\ref{1to2}) of \S\ref{remain} applied to all finite
subextensions $k'|k$ of $K|k$ yields a surjection
\[
\ressum_{\p \in S_\R \smallsetminus S (K)}
H^1(T(K_\p(2)|K_\p),\Q_2/\Z_2) \surjr{}
H^2(G(K_S(2)|K),\Q_2/\Z_2).
\]
and therefore, in order to prove the proposition, it suffices to
show that the group $H^2(G(K_S(2)|K),\Q_2/\Z_2)$ is $2$-divisible.
This is trivial if $S \cap S_\R(K)=\varnothing$ because then
$\cd\, G(K_S(2)|K)\leq 2$. Otherwise, this follows from the
commutative diagram
\[
\renewcommand{\arraystretch}{1.5}
\begin{array}{ccc}
H^2(G(K_S(2)|K),\Q_2/\Z_2)/2 & \hookrightarrow &\!\!\!H^3(G(K_S(2)|K),\Z/2\Z)\\
\mapd{} &&\mapd{\wr}\\
\ressum_{\p \in S \cap S_\R (K)}\!\!\!\!
H^2(T(K_\p(2)|K_\p),\Q_2/\Z_2)/2& \hookrightarrow
&\!\!\!\ressum_{\p \in S \cap S_\R (K)} \!\!\!\!
H^3(T(K_\p(2)|K_\p),\Z/2\Z).
\end{array}
\renewcommand{\arraystretch}{1.5}
\]
The right hand vertical arrow is an isomorphism by
proposition~\ref{isoab3}. But $H^2(T(K_\p(2)|K_\p),\Q_2/\Z_2)=0$
for all $\p\in S \cap S_\R (K)$ and therefore the object in the
lower left corner is zero.
\enddemo

\section{Relation to \'{e}tale cohomology} \label{4}

Let $k$ be a number field and $S$ a finite set of places of $k$. We think of $\Spec( \Cal
O_{k,S})$ as  ``$\{\hbox{scheme-theoretic points of $\Spec( \Cal O_{k,S})$}\} \cup
\{\hbox{real places of $k$ not in $S$\}}$''.\linebreak Essentially following Zink \cite{Z},
we introduce the site $\Spec ( \Cal O_{k,S})_{\et,\md}$.

\smallskip
{\it Objects} of the category are pairs $\bar U=(U,U_\real)$, where $U$ is a scheme together
with an \'{e}tale structural morphism $\phi_U: U \to \Spec(\Cal O_{k,S})$ and  $U_\real$ is a
subset of the set of real valued points $U(\R)=\Mor_{\Schemes}(\Spec(\R),U)$ of $U$ such
that $\phi_U(U_\real) \subset S_\R(k) \smallsetminus S$.

\smallskip
{\it Morphisms} are scheme morphisms $f: U \to V$ over $\Spec(\Cal O_{k,S})$ satisfying
$f(U_\real)\subset V_\real$.

\smallskip
{\it Coverings} are families $\{\pi_i: \bar U_i \to \bar U\}_{i\in I}$ such that $\{\pi_i:
U_i \to U\}_{i\in I}$ is an \'{e}tale covering in the usual sense and $\bigcup_{i\in I}
\pi_i({U_i}_\real)=U_\real$.

\medskip\noindent
There exists an obvious morphism of sites
\[
\Spec( \Cal O_{k,S})_{\et} \lang \Spec( \Cal O_{k,S})_{\et,\md}
\]
and both sites coincide if $S$ contains all real places of $k$.
 The pair $\bar X=(X,X_\real)$ with  $X=\Spec(\Cal
O_{k, S})$ and $X_\real=S_\R(k)\smallsetminus S$ is the terminal object of the category and
the profinite group $G_S(k)$ is nothing else but the fundamental group of $\bar X$ with
respect to this site. Let $\eta$ denote the generic point of $X$. For a sheaf $A$ of abelian
groups on $\Spec( \Cal O_{k,S})_{\et,\md}$ and for any point $v$ of $\bar X$ we have a
specialization homomorphisms $s_v: A_v \to A_\eta$ from the stalk  $A_v$ of $A$ in $v$ to
that in $\eta$. For each point $v\in X_\real$ we consider the local cohomology $H_v^i(\bar
X, A)$ with support in $v$. There is a long exact localization sequence (see \cite{Z})
\[
\cdots \to \bigoplus_{v \in X_\real} H_v^i(\bar X, A) \to H^i_{\et,\md}(\bar X,A) \to
H^i_\et(X,A) \to \cdots
\]
and the local cohomology with support in  real points is calculated as follows:

\begin{lemma}\label{loccoh} For $v \in X_\real$ the following holds.
\[
\renewcommand{\arraystretch}{1.2}
\begin{array}{ccll}
H^0_v(\bar X, A)&=& \ker (s_v: A_v \to A_\eta)&\\
H^1_v(\bar X, A)&=& \coker (s_v: A_v \to A_\eta)&\\
H^i_v(\bar X, A)&=& H^{i-1}(k_v, A_v)& \hbox{ for } i \geq 2.
\end{array}
\renewcommand{\arraystretch}{1}
\]
Here the right hand side of the last isomorphism is the Galois cohomology of the field $k_v$.
\end{lemma}

\demo{Proof} See \cite{Z}, Lemma 2.3. \enddemo

\noindent {\bf Remark:} Suppose that $S$ contains all primes dividing $2$ and no real primes.
Let $A$ be a locally constant constructible sheaf on $\Spec(\Cal O_{k,S})_\et$ which is
annihilated by a power of $2$. We denote the push-forward of $A$ to $\Spec(\Cal
O_{k,S})_{\et,\md}$ by the same letter. By Poitou-Tate duality, the boundary map of the long
exact localization sequence
\[
H^i_\et(X,A) \lang \bigoplus_{v \in X_\real} H_v^{i+1}(\bar X, A)= \bigoplus_{v \hbox{ \rm
\scriptsize arch.}} H^{i}(k_v, A_v)
\]
is an isomorphisms for $i \geq 3$ and surjective for $i=2$. Therefore, we obtain the
vanishing of $H^i_{\et,\md}(\Spec(\Cal O_{k,S}),A)$ for $i \geq 3$. In this situation  the
modified \'{e}tale cohomology is connected to the ``positive \'{e}tale cohomology''
$H^*_{2}(\Spec(\Cal O_{k,S}),A_+)$ defined in \cite{C} in the following way. There exists a
natural exact sequence

\medskip\noindent $ 0  \to H^0_{\et,\md}(\Spec(\Cal
O_{k,S}),A) \to$
\[
\bigoplus_{v\, \hbox{\rm\scriptsize arch.}} H^0(k_v,A_v) \to H^0_{2}(\Spec(\Cal O_{k,S}),A_+)
\to H^1_{\et,\md}(\Spec(\Cal O_{k,S}),A) \to 0.
\]
and isomorphisms
\[
H^{i}_{2}(\Spec(\Cal O_{k,S}),A_+)\mapr{\sim} H^{i+1}_{\et,\md}(\Spec(\Cal O_{k,S}),A)
\]
for $i \geq 1$. This can be easily deduced from the long exact localization sequence, lemma
\ref{loccoh} and  the long exact sequence (2.4) of \cite{C}.

\bigskip
Now let $A$ be a discrete $G_S(k)$-module. The module $A$ induces locally constant sheaves on
$\Spec (\Cal O_{k,S})_{\et, \md}$ and $\Spec (\Cal O_{k,S})_{\et}$, which we will denote by
the same letter. According to lemma \ref{loccoh}, we obtain for every $v \in X_\real$
\[
H^i_v(\bar X,A)=0 \quad \hbox{ for } i=0,1.
\]
Let $\widetilde X = (\Spec (\Cal O_{k_S,S}), S_\R(k_S)\smallsetminus S(k_S)$ be the universal
covering of $\bar X$. The Hochschild-Serre spectral sequence
\[
 E_2^{ij}= H^i(G_S(k),H^j_{\et,\md}(\widetilde X,A) ) \Longrightarrow
H^{i+j}_{\et,\md}(\bar X,A)
\]
induces natural comparison homomorphisms
\[
H^i(G_S(k),A) \lang H^i_{\et,\md}(\bar X,A)
\]
for all $i \geq 0$. It follows immediately from the spectral sequence that these
homomorphisms are isomorphisms if
\[
H^j_{\et,\md}(\widetilde X,A)= 0
\]
for all $j \geq 1$.

\bigskip Next we are going to prove theorem \ref{et} of the introduction. Assume that $S$
contains all primes dividing $2$ and that $A$ is $2$-torsion. Both sides of the comparison
homomorphism commute with direct limits, and so, in order to prove theorem~\ref{et}, we may
suppose that $A$ is finite. Since $A$ is constant on $\widetilde X$, we can easily reduce to
the case $A=\Z/2\Z$, in order to show $H^j_{\et,\md}(\widetilde X,A)=0$ for $j \geq 1$.
Furthermore, the assertion is trivial for $j=1$. The theorem is well-known if $S$ contains
all real primes (see \cite{Z}, prop. 3.3.1 or \cite{M}, II, 2.9) and so, passing to the
limit over all finite subextensions of $k$ in $k_S$, we obtain natural isomorphisms for all
$j \geq 0$.
\[
H^j(G_{S\cup S_\R}(k_S),\Z/2\Z) \mapr{\sim} H^j_\et(\widetilde X \smallsetminus
S_\R(k_S),\Z/2\Z).
\]
On the other hand, theorem~\ref{riemann} for $T=S\cup S_\R$, $S=S$ applied to all finite
subextensions  of $k$ in $k_S$ in conjunction with theorem~\ref{prof} induces isomorphisms
for all $j\geq 1$.
\[
H^j(G_{S\cup S_\R}(k_S),\Z/2\Z) \mapr{\sim} \ressum_{v \in S_\R\smallsetminus S(k_S)}
H^j(k_v,\Z/2\Z).
\]
These two isomorphisms together with the long exact localization sequence show that
\[
H^j_{\et,\md}(\widetilde X, \Z/2\Z)=0
\]
for $j\geq 1$. This completes the proof of theorem \ref{et}.

\bigskip\noindent
Theorem~\ref{et} is best understood in the context of \'{e}tale homotopy, namely as  a vanishing
statement on the $2$-parts of higher homotopy groups. For a scheme $X$ we denote by $X_{et}$
its \'{e}tale homotopy type, i.e.\ a pro-simplicial set. The \'{e}tale homotopy groups of $X$ are by
definition the homotopy groups of $X_{et}$ and, as is well known, these pro-groups are
pro-finite, whenever the scheme $X$ is noetherian, connected and geometrically unibranch
(\cite{artinmazur} Theorem 11.1). If we consider the modified \'{e}tale site $\Spec ( \Cal
O_{k,S})_{\et,\md}$ as above, we obtain in exactly the same manner as for the usual \'{e}tale
site a pro-finite simplicial set $\bar X_{\et,\md}$. We denote the universal covering of
$\bar X_{\et,\md}$ by $\widetilde{X}_{\et,\md}$. If $p$ is a prime number and $Y_\cdot$ is a
pro-simplicial set, we denote the pro-$p$ completion of $Y_\cdot$ by $Y_\cdot^{\wedge p}$.
Furthermore, we write $G(p)$ for the  maximal pro-$p$ factor group of a pro-group $G$.

\begin{lemma}\label{pro-p-ver}
Assume that $Y_\cdot$ is a simply connected (i.e.\ $\pi_1(Y_\cdot)=0$) pro-sim\-pli\-cial set
such that $\pi_i(Y_\cdot)$ is pro-finite for all $i\geq 2$. Then we have isomorphisms for
all $i$:
\[
\pi_i(Y_\cdot )(p) \mapr{} \pi_i(Y_\cdot^{\wedge p}).
\]
\end{lemma}

\demo{Proof:} See \cite{S}, prop.\ 13. \enddemo

\noindent
 For a pro-group $G$ we denote by $K(G,1)$ the Eilenberg-MacLane pro-simplicial set
associated with $G$ (cf.\ \cite{artinmazur}, (2.6)). If $S$ contains all real primes of $k$
the following theorem was proved in \cite{S}, prop.\ 14.

\begin{theorem} Let $k$ be a number field and $S$ a finite set of primes of\/ $k$ containing all primes
dividing $2$. Let $\bar X$ be the pair $(X,X_\real)$ with  $X=\Spec(\Cal O_{k, S})$ and
$X_\real=S_\R(k)\smallsetminus S$ endowed with the modified \'{e}tale topology. Then the higher
homotopy groups of\/ $\bar X_{\et,\md}$ have no $2$-part, i.e.
\[
\pi_i(\bar X_{\et,\md})(2)=0 \quad \hbox{ for } i \geq 2.
\]
Furthermore,  the canonical morphism
\[
(\bar X_{\et,\md})^{\wedge 2} \mapr{} K(G_S(k)(2),1)
\]
is a weak homotopy equivalence.
\end{theorem}

\demo{Proof:}  Since $G_S(k)$ is the fundamental group of $\bar X_{\et,\md}$,
theorem~\ref{et} implies that the universal covering $\widetilde{X}_{\et,\md}$ of $\bar
{X}_{\et,\md}$ has no cohomology with values in $2$-primary coefficient groups. By the
Hurewicz theorem (\cite{artinmazur}, (4.5)),  the pro-$2$ completion of $\widetilde
X_{\et,\md}$ is weakly contractible. Therefore lemma \ref{pro-p-ver} implies
 $$\pi_i(\bar X_{\et,\md})(2)\cong \pi_i(\widetilde
X_{\et,\md})(2)\cong \pi_i\big((\widetilde X_{\et,\md})^{\wedge 2}\big)=0$$
 for $i \geq 2$, which shows
the first statement of the theorem. By theorem~\ref{prof}, for every finite $2$-primary
$G_S(k)(2)$-torsion module $A$ the inflation homomorphism $H^{i}(G_S(k)(2),A) \mapr{}
H^{i}(G_S(k),A)$ is an isomorphism for all $i$. The same arguments as above show that the
universal covering of $(\bar X_{\et,\md})^{\wedge 2}$ is weakly contractible. This proves the
second statement.
\enddemo

\section{Closing Remarks}
\underline{1. {\em Dualizing modules}}

\medskip\noindent
Unfortunately, we do not have (despite semi-tautological
reformulations of the definition) a good description of the
$p$-dualizing module $I$ of the group $G_{S}$, where $S$ is a
finite set of finite primes containing $S_p$. If  $k$ is totally
imaginary, then $I$ is determined by the exact sequence
\[
0 \lang \mu_{p^\infty}  \mapr{\hbox{\scriptsize\it diag}}
\ressum_{\p \in S(k_{S})} \mu_{p^\infty} \lang I \lang 0
\]
(see \cite{NSW}, (10.2.1)) and the group $G_S$ is a duality group
at $p$ of dimension $2$ (see \cite{S}, th.4 or \cite{NSW},
(10.9.1)). The general case remains unsolved (also for odd $p$).

\bigskip\noindent
\underline{2. {\em Free profinite product decompositions}}

\medskip\noindent
In this paper we used free pro-$p$-product decompositions of
Galois groups of pro-$p$-extensions of global fields into  Galois
groups of local pro-$p$-extensions in an essential way. One might
ask whether, for sets of places $T \supset S$, the natural
homomorphism
\[
\phi: \freeproductmed_{\p \in T \smallsetminus S (k_S)} T(\bar
k_\p|k) \lang G(k_T|k_S)
\]
is an isomorphism, where the free product on the left hand side is
the free product of {\it profinite} groups. More precisely, one
has to ask, whether there exists a continuous section to the
natural projection $T\smallsetminus S (k_T) \to T \smallsetminus
S (k_S)$ such that the above map is an isomorphism (cf.\ the
discussion in section~\ref{secfree}). We do not know the answers
to this question in general. It is \lq yes\rq\ if $S$ contains all
but finitely many primes of $k$ (see below). But it seems likely
that $\phi$ is never an isomorphism if $T$ and $S$ are finite.
The present level of knowledge on this question is rather low.
For example, we do not know whether there are infinitely many
prime numbers $p$ such that $p^\infty$ divides the order of
$G_T$. The best result known in this direction is that if $T$
contains all real places and all primes dividing one prime number
$p$, then there exist infinitely many prime numbers $\ell$
dividing the order of $G_T$ (see \cite{Wa}, cor.3 or \cite{NSW},
(10.9.4)).

In the case that $S$ contains all but finitely many primes of
$k$, we can deduce the above statement by applying the following
slightly more general result to the complement of $S$:

For a finite set $S$ of primes of $k$, let $k^S$ be the maximal
extension of $k$ in which all primes in $S$ are totally
decomposed. Then there exists a continuous section to the natural
projection $S(\bar k)\to S(k^S)$ such that the natural map
\[
\freeproductmed_{\p \in S(k^S)} G(\bar k_\p|k) \lang G(\bar k
|k^S)
\]
is an isomorphism. This had been proved first in the special case
$S=S_\R$ by Fried-Haran-V\"{o}lklein \cite{FHV} and then by Pop
\cite{P} for arbitrary finite $S$.

\bigskip\pagebreak[3]

\noindent
\underline{3. {\em Leopoldt's conjecture}}

\medskip\noindent
The Leopoldt conjecture for $k$ and a prime number $p$ holds if
and only if the group
\[
H^2(G_S,\Q_p/\Z_p)
\]
is trivial for one (all) finite set(s) of primes $S \supseteq
S_p$. The weak Leopoldt conjecture is true for $k$, $p$ and a
$\Z_p$-extension $k_\infty|k$ if and only if
\[
H^2(G_S(k_\infty),\Q_p/\Z_p)
\]
is trivial for one (all) finite set(s) of primes $S \supseteq S_p$
(of $k$). This is well known for odd $p$ and for $p=2$ it can be
easily deduced from the above results.

\bigskip \pagebreak[3]

\noindent \underline{4. {\em Iwasawa theory}}

\medskip\noindent
Let $k$ be a number field, $S\supseteq S_2$ a finite set of
primes of $k$ and $k_\infty|k$ the cyclotomic $\Z_2$-extension of
$k$. Let $\varGamma=G(k_\infty|k)\cong \Z_2$ and let $\varLambda=
\Z_2\wingl \varGamma \wingr \cong \Z_2[[T]]$ be the Iwasawa
algebra. We consider the compact $\varLambda$-module
\[
X_S= G(k_S(2)|k_\infty)^{ab}.
\]
 Then the
following holds\medskip

\Item{\rm (i)} $X_S$ is a finitely generated
$\varLambda$-module.\smallskip \Item{\rm (ii)} $\rank_\varLambda
\, X_S = r_2$ (the number of complex places of $k$).
\smallskip \Item{\rm (iii)} $X_S$ does not contain any nontrivial
finite $\varLambda$-submodule. \smallskip \Item{\rm (iv)} the
$\mu$-invariant of $X_S$ is greater than or equal to $\#\; S \cap
S_\R (k)$.

\medskip\noindent
Properties (i)-(iii) follow in a purely formal way  (see
\cite{NSW}, (5.6.15)) from the facts that:  (a)
$\chi_2(G_S(2))=-r_2$, (b) $H^2(G_S(k_\infty)(2),\Q_2/\Z_2)=0$ and
(c) $H^2(G_S(2),\Q_2/\Z_2)$ is $2$-divisible. Assertion (iv) is
trivial if  $S$ contains no real places and in the general case
it follows from the exact sequence
\[
0 \to (\varLambda/2)^{\#\, S\cap S_\R(k)} \lang X_S  \lang X_{S
\smallsetminus S_\R} \lang 0.
\]

\medskip\noindent
Now let $k^+$ be a totally real number field, $k=k^+(\mu_{4})$,
$k^+_\infty$ the cyclotomic $\Z_2$-extension of $k^+$ and
$k_\infty=k^+_\infty(\mu_{4})=k(\mu_{2^\infty})$. Let $k_n$ be
the unique subextension of degree $2^n$ in $k_\infty$ and let $J$
be the complex conjugation. We set $A_n= \Cl(k_n) (2)$ and
\[
A_n^- := \{a \in A_n\, |\, a J(a)=1 \}.
\]
Furthermore, let $A_\infty^-= \dirlim A_n^-$, $X^+=X_{S_2}(k^+)$,
let $\scriptstyle \vee$ denote the Pontryagin dual and  $(-1)$
the Tate-twist by $-1$. Then there exists a natural homomorphism
\[
\phi: (A_\infty^-)^\vee \lang X^+(-1)
\]
whose kernel and cokernel are annihilated by $2$. If the Iwasawa
$\mu$-invariant of $k$ is zero (this is known if $k|\Q$ is
abelian), then $\phi$ is a pseudo-isomorphism, i.e.\ $\phi$ has
finite kernel and cokernel. This can be seen by a slight
modification of the arguments given in \cite{G}, \S2:

Let $M^+$ be the maximal abelian $2$-extension of $k^+_\infty$
which is unramified outside $S_2$, in particular, $M^+$ is totally
real. Kummer theory shows that, for an $\alpha \in
k_\infty^\times$, the field $k_\infty(\sqrt[2^n]{\alpha})$ is
contained in $M^+k_\infty$ if and only if: (a) $\alpha \in
k_{\infty,\p}^{\times 2^n}$ for all $\p \not \in S(k_\infty)$ and
(b) $\alpha J(\alpha)= \beta^{2^n}$ for a totally positive element
$\beta \in k_\infty^+$. Let $R_n$ be the subgroup in
$k^\times_\infty/k^{\times 2^n}_\infty$ generated by elements
satisfying (a) and (b) and let
$$
\mathfrak M^- := \dirlim_n R_n \subset k_\infty^\times \,
\otimes\, \Q_2/\Z_2.
$$
Then we have a perfect Kummer pairing $X^+ \times {\mathfrak M}^-
\to \mu_{2^\infty}$. Since all primes dividing $2$ are infinitely
ramified in $k_\infty|k$, for $\alpha \otimes 2^{-n} \in
{\mathfrak M}^-$ there exists a unique ideal $\mf a$ in
$k_\infty$ with $\mf a^{2^n} =(\alpha)$ and the class $[\,\mf
a\,]$ is contained in $A_\infty^-$. This yields a homomorphism
\[
\phi^\vee: \mf M^- \lang A_\infty^-.
\]
A straightforward computation shows that $\im(\phi^\vee)\supseteq
(A_\infty^-)^2$ and that $\ker (\phi^\vee)$ is the image of $\Cal
O^\times_{k_\infty^+,\varnothing}/\Cal O^{\times
2}_{k_\infty^+,S_\R}$ in $\mf M^-$ (notational conventions as in
\S\ref{remain}). Thus, if the Iwasawa $\mu$-invariant of $k$ is
zero, then the cokernel of $\phi^\vee$ is finite and it remains to
show the same for its kernel. Since $\mu=0$, the $\F_2$-ranks of
$\null_2\Cl^0(k_n^+)$ (the subgroup of elements annihilated by
$2$ in the ideal class groups in the narrow sense) are bounded
independently of $n$. Thus also the $\F_2$-ranks of the kernels
of the signature maps
\[
\Cal O^\times_{k_n^+,S_\R}/\Cal O^{\times 2}_{k_n^+,S_\R} \lang
\bigoplus_{v\in S_\R(k_n^+)} \R^\times/ \R^{\times 2}
\]
are bounded independently of $n$. But the direct limit over $n$
of these kernels is just the group in question. Finally, we
obtain the result by taking Pontryagin duals.

\vfill

\noindent {\sc  Alexander Schmidt}
  Mathematisches Institut,
  Universit\"{a}t Heidelberg,
  Im Neuenheimer Feld 288,
  69120 Heidelberg,
  Deutschland\\
  e-mail: schmidt@mathi.uni-heidelberg.de

\end{document}